\documentclass[twocolumn]{svjour3}  
\smartqed  

\usepackage[dvipdfmx]{color}
\usepackage{graphicx}
\usepackage[T1]{fontenc}
\usepackage{mathptmx}  
\usepackage[scaled]{helvet}  
\usepackage{courier}
\usepackage[subrefformat=parens]{subcaption}

\makeatletter
\let\cl@chapter\undefined
\makeatletter
\usepackage{amsmath,amssymb}   
\usepackage{mathtools} 
\usepackage{cleveref}  
\Crefname{equation}{Eq.}{Eqs.}%
\Crefname{figure}{Fig.}{Figs.}%
\usepackage{amsbsy}
\usepackage{bm}
\usepackage{algorithmic}
\usepackage{algorithm}

\begin{document}

\title{Approximate probability density function for nonlinear surging in irregular following seas
%
}

\author{Atsuo Maki         \and Yuuki Maruyama \and
        Keiji Katsumura         \and Leo Dostal
}

\institute{Atsuo Maki \and Yuuki Maruyama \and Keiji Katsumura \at
              Osaka University, 2-1 Yamadaoka, Suita, Osaka, Japan \\
              \email{maki@naoe.eng.osaka-u.ac.jp} 
           \and
           Leo Dostal \at
           Institute of Mechanics and Ocean Engineering, Hamburg University of Technology, 21043 Hamburg, Germany\\
           \email{dostal@tuhh.de} 
}

\date{Received: date / Accepted: date}

\maketitle

\begin{abstract}
The broaching that follows the surf-riding is a dangerous phenomenon that can lead to the capsizing of a vessel due to its violent yaw motion. Most of the previous studies on surf-riding phenomena in irregular waves have been conducted by replacing irregular waves with regular waves. In contrast, this study provides suggestions on how to directly calculate nonlinear surge motion in irregular seas. In this study, the statistical aspects of the surf-riding phenomenon are first presented. Then, under several approximations, we show how to calculate the probability density function theoretically. Although the results obtained are based on strong approximations, it is found that the nonlinear surge oscillations in irregular following seas can be explained from a qualitative point of view.

\keywords{Irregular seas \and Stochastic differential equation \and Broaching-to \and Bifurcation}
\end{abstract}

\section{Introduction}
    A ship operating in following/quartering seas may occasionally experience surf-riding and subsequent broaching. The dangers of this phenomenon are well known\cite{saunders1965}, and various experimental, theoretical, and numerical approaches have been conducted \cite{Grim1951,Makov1969,Ananiev1966,Kan_1990_surfriding,umeda1990SR_regular,umeda1992Broaching,spyrou1996dynamic,spyrou1996dynamic2,spyrou1997dynamic,umeda1999nonlinear}. 
    
    From the aspect of a nonlinear dynamical system, the surf-riding boundary corresponds to a heteroclinic orbit. Unfortunately, the nonlinear surge equation cannot be solved analytically in general. Therefore, for a theoretical approach, it is necessary to introduce approximate methods such as Melnikov's method~\cite{holmes1980averaging}. Kan applied Melnikov's method to propose an approximate formula for predicting the surf-riding threshold~\cite{Kan_1990_surfriding}. Spyrou~\cite{spyrou2006} extended Kan's approach and obtained an approximate formula for the surf-riding threshold. Further, Maki generalized their approaches and obtained an approximate formula~\cite{maki2010surfriding}. A comprehensive review on theoretical estimation methods for surf-riding thresholds in regular waves is available in the literature \cite{maki2024NOLTA}.

    In regular seas, the wave celerity can be rigorously determined, and the surf-riding can be defined by the phenomenon in which the ship travels at the same speed as the wave celerity. In irregular seas, on the other hand, the wave celerity changes constantly with time and location, making it difficult to define the wave celerity of waves in a practical sense. Therefore, in Umeda's work \cite{umeda1990SR_regular}, the probability of surf-riding occurrence is determined by approximately replacing irregular waves with regular waves. Maki and Miyauchi \cite{Maki2015_FAST} conducted free-running model experiments in a 247m-long towing tank at the Naval Systems Research Center of Acquisition, Technology \& Logistics Agency and then attempted to define the surf-riding phenomenon in irregular waves based on the trend of the pitch angle of the ship model. Belenky et al. \cite{belenky2012,belenky2016SR} examined the reduction of overturning due to broaching and wave riding using a split-time method approach. Dostal and Kreuzer~\cite{dostal2013SR} used Melnikov's method\cite{simiu2002chaotic} for irregular disturbances and conducted a theoretical study on the limits of surf-riding in irregular seas.
    
    On the other hand, Spyrou et al. \cite{spyrou2014detection} were actively studying the ``Instantaneous wave celerity'' and used it to approach the surf-riding phenomenon from that perspective. This idea emerges from the use of the ``Instanteneosu frequency''. A corresponding comprehensive review can be found in \cite{boashash1992estimating}. Moreover, Spyrou et al. \cite{kontolefas2020probability,spyrou2023SR} have introduced the concept of ``high-run''. This is the percentage of time that the ship's speed is above a threshold value. In \cite{spyrou2023SR}, also the instantaneous wave celerity is used. Then, probability estimates of high-run are calculated. In the research group of Spyrou, research on surf-riding phenomena in irregular waves from the viewpoint of wave celerity or using approaches based on nonlinear dynamical systems have been conducted, as can be seen in the references \cite{spyrou2012conditions,spyrou2014surf,spyrou2018nonlinear,tigkas2023hybrid}.
    
    However, all of them were either experiments only, fully/partly numerical calculations, approaches that replaced irregular waves with regular waves, or approaches in which the wave celerity in irregular seas is alternatively used. So far, we have not found studies that discuss the issue from a fully probabilistic perspective. Therefore, in this study  the nonlinear surge motion of a ship in irregular waves is investigated from a statistical point of view. After this a theoretical approach based on the theory of stochastic processes is presented.

    Initial results from the investigation described in this study were developed by Maki et al. \cite{Maki2023SR_JASNAOE}. In this paper, the results are presented more extensively, with more details and some revisions.

\section{Notations}\label{sec:notaions}
    In this study, the $n$-dimensional Euclidean space is denoted by $\mathbb{R}^n$, and the set of real numbers for $n=1$ is denoted by $\mathbb{R}$. The expectation operation is denoted by $\mathbb{E}$, and $t$ represents time. The overdot of time-dependent variables indicates the derivative with respect to time $t$.

\section{Subject ship}\label{sec:subject_ship}
    The subject ship in this study is the DTMB5415 hull form (ship length: $L_{\mathrm{PP}}=2.75~\mathrm{m}$), and it is the same subject ship used in the literature \cite{maki2010surfriding,maki2024NOLTA}. 
    \begin{figure}[htb]
        \centering 
        \includegraphics[width=1.0\hsize]{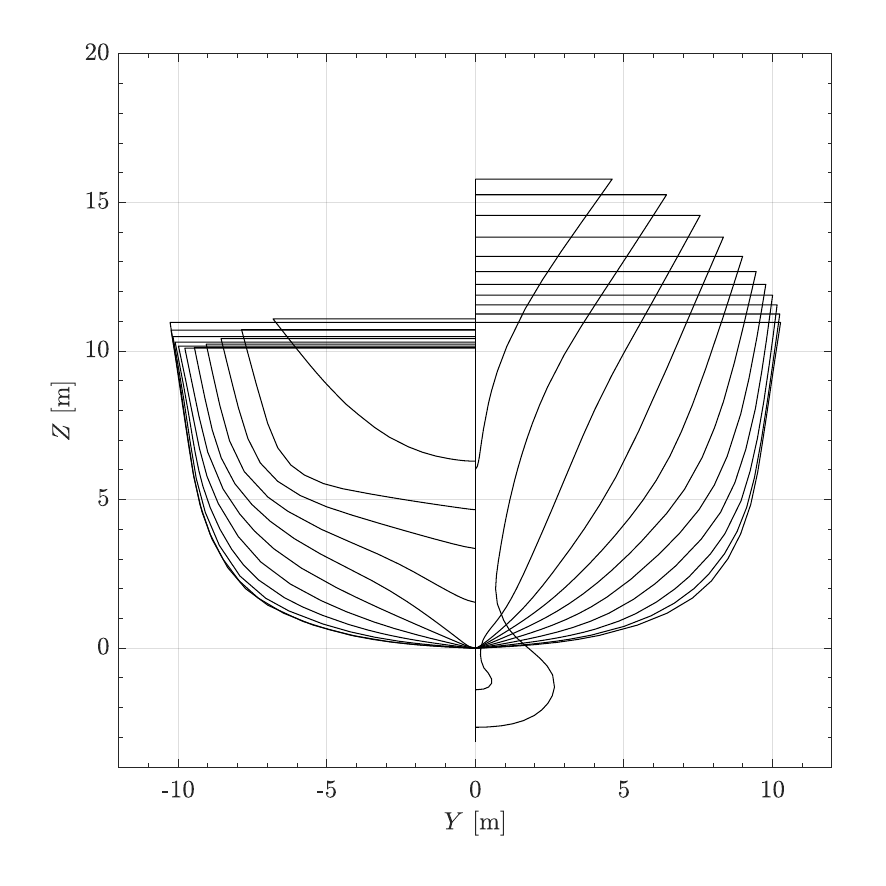}
        \caption{Body plan of the subject ship (DTMB5415). Here, ship length is $L_{\mathrm{PP}}=2.75~\mathrm{m}$. The detailed description of this ship can be found in the literature \cite{maki2016surfriding,maki2024NOLTA}.}
        \label{fig:DTMB5415}
    \end{figure}
    It is well-known that there exist two types of thresholds for surf-riding. One is called the ``surf-riding threshold'' which exists in lower-speed regions. The other one is called the ``wave blocking threshold'' which exists in higher speed regions \cite{maki2014surfriding}. Here, the Froude numbers, which correspond to these two thresholds, are represented as $\overline{\mathrm{Fn}}_{\mathrm{LWR}}$ and $\overline{\mathrm{Fn}}_{\mathrm{UPS}}$, respectively. The Froude number which corresponds to the ship speed in calm water with the same propeller revolutions is defined by $\overline{\mathrm{Fn}}$. Here, the authors define the ship speed in calm water as $\bar{u}_{\mathrm{S}}$. Then, $\overline{\mathrm{Fn}}$ can be calculated as:
    \begin{equation}
        \mathrm{\overline{Fn}}=\frac{\bar{u}_{\mathrm{S}}}{\sqrt{L_{\mathrm{S}}g}}  
        \label{eq:def_Fn}
    \end{equation}
    where $g$ is the gravitational acceleration and $L_{\mathrm{S}}$ is the ship length.
    
\begin{figure}[htb]
        \centering 
        \includegraphics[width=1.0\hsize]{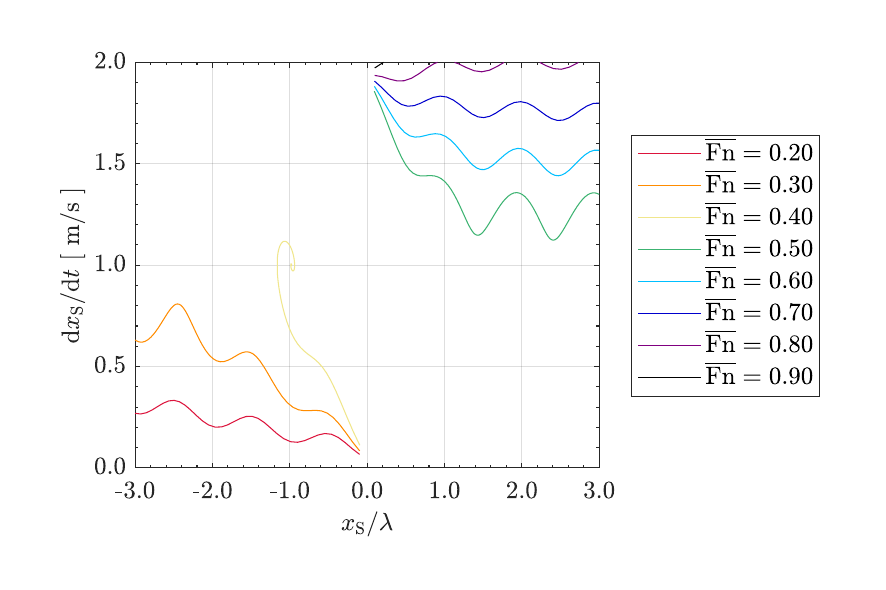}
        \caption{Phase portraits for DTMB5415 with $H/\lambda =0.04$ and $\lambda/L_{\mathrm{S}} =1.0$}
        \label{fig:phase_portraits_1p0}
    \end{figure}

\begin{figure}[htb]
        \centering 
        \includegraphics[width=1.0\hsize]{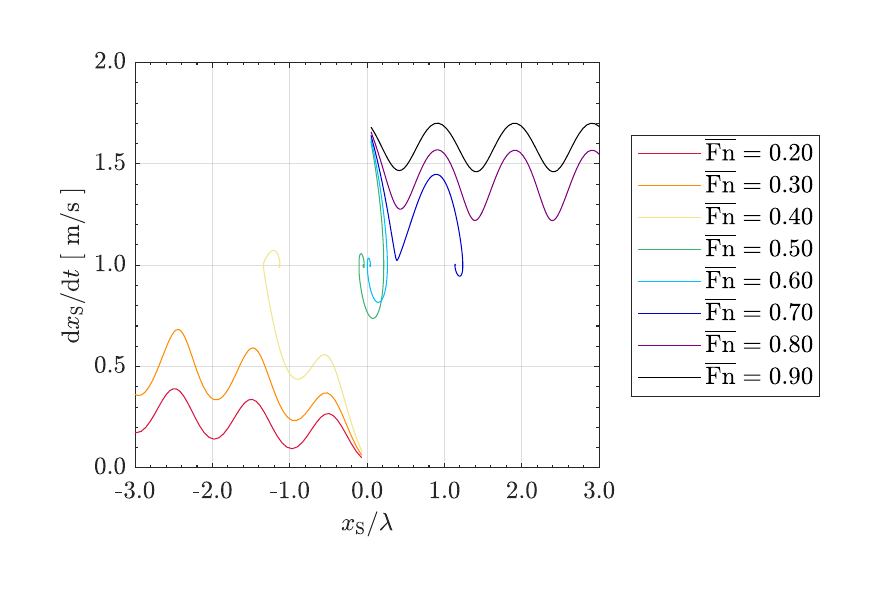}
        \caption{Phase portraits for DTMB5415 with $H/\lambda =0.04$ and $\lambda/L_{\mathrm{S}} =2.0$}
        \label{fig:phase_portraits_2p0}
    \end{figure}

    The surf-riding thresholds can be theoretically estimated By utilizing the method provided by Maki et al.\cite{maki2010surfriding}. The estimation formula is based on Melnikov's method and the formula provided in the literatue \cite{maki2010surfriding} was adopted to second generation intact stability criteria by the International Maritime Organization (IMO). The estimated surf-riding thresholds are shown in \Cref{tab:SR_threshold}. Here, $\lambda$ is the wavelength of the incident wave, and $H$ is the wave height.
    The phase portraits for different $\overline{\mathrm{Fn}}$ are shown in \Cref{fig:phase_portraits_1p0,fig:phase_portraits_2p0}. \Cref{fig:phase_portraits_1p0} and \Cref{fig:phase_portraits_2p0} show the results for $\lambda/L_{\mathrm{S}}=1.0$ and $\lambda/L_{\mathrm{S}}=2.0$, respectively. In both figures, there are two cases where the ship periodically overtakes or is overtaken by waves ($\overline{\mathrm{Fn}}=0.2,~0.3,~0.8,~0.9$ in Fig.~\ref{fig:phase_portraits_2p0}). These periodic cases are denoted as \textbf{Case 1}. On the other hand, the ship motion finally converges to the equilibrium point  ($\overline{\mathrm{Fn}}=0.4,~0.5,~0.6,~0.7$ in Fig.~\ref{fig:phase_portraits_2p0}). These periodic cases are denoted as \textbf{Case 2}. Here, its thresholds are the surf-riding and wave-blocking thresholds. This is due to a global qualitative change in the dynamical system caused by the change in $\overline{\mathrm{Fn}}$. Therefore, in the framework of nonlinear dynamical system theory, this threshold is called a bifurcation. The bifurcation in the present case is a heteroclinic bifurcation. 
    \begin{table}[h]
        \caption{Surf-riding threshold for $H/\lambda = 0.04$}
        \centering
        \begin{tabular}{lc|lc} 
        $\lambda/L_{\mathrm{P}}$ & $\overline{\mathrm{Fn}}_{\mathrm{LWR}}$ & $\lambda/L_{\mathrm{P}}$ & $\overline{\mathrm{Fn}}_{\mathrm{UPS}}$ \\
        \hline 
        $1.0$ & $0.3166$ & $1.0$ & $0.4756$ \\
        $2.0$ & $0.3900$ & $2.0$ & $0.7305$ 
        \end{tabular}
        \label{tab:SR_threshold}
    \end{table}

\section{Coordinate systems and equation of motion}\label{sec:Equation_of_Motion}
    The ship is assumed to perform only surge motion and the coordinate of the center of gravity in the earth-fixed coordinate system is $x_{\mathrm{S}}$. Then, the equation of the surge motion of the vessel is represented as:
    \begin{equation}
        \label{eq:EoM_original}
        (m+m_x)\frac{\mathrm{d}^2x_{\mathrm{S}}}{\mathrm{d}t^2} = - R\left(\frac{\mathrm{d}x_{\mathrm{S}}}{\mathrm{d}t}\right) + T\left(\frac{\mathrm{d}x_{\mathrm{S}}}{\mathrm{d}t}, n_{\mathrm{P}}\right) + f(x_{\mathrm{S}}, t)
    \end{equation}
    Here, $m$: the ship mass, $m_x$: the added mass in the surge direction, $n_{\mathrm{P}}$: the propeller revolution number, $R(\mathrm{d}x/\mathrm{d}t)$: the ship resistance, $T(\mathrm{d}x/\mathrm{d}t, n_{\mathrm{P}})$: the propeller thrust, $f(x_{\mathrm{S}}, t)$: wave induced surge force. It is assumed that the ship resistance $R(\mathrm{d}x/\mathrm{d}t)$ can be approximated by a 5th order polynomial of $\mathrm{d}x/\mathrm{d}t$. 
    \begin{equation}
        \label{eq:polynomial_resistance}
        R\left(\frac{\mathrm{d}x_{\mathrm{S}}}{\mathrm{d}t}\right) \approx \sum_{i=1}^{5} r_i \left(\frac{\mathrm{d}x_{\mathrm{S}}}{\mathrm{d}t}\right)^i
    \end{equation}
    The propeller thrust $T(\mathrm{d}x/\mathrm{d}t, n_{\mathrm{P}})$ can be calculated by using $K_{\mathrm{T}}$. 
    \begin{equation}
        T\left(\frac{\mathrm{d}x}{\mathrm{d}t}, n_{\mathrm{P}}\right) = (1-t_{\mathrm{P}})\rho n_{\mathrm{P}}^2 D_{\mathrm{P}}^4 K_{\mathrm{T}}(J_{\mathrm{T}}) 
    \end{equation}
    Here, $t_{\mathrm{P}}$: the thrust deduction coefficient, $\rho$: water density, $D_{\mathrm{P}}$: the propeller diameter. $J_{\mathrm{T}}$ can be defined as:
    \begin{equation}
        J_{\mathrm{T}} = \frac{(1-w_{\mathrm{P}})}{n_{\mathrm{P}} D_{\mathrm{P}}}\frac{\mathrm{d}x_{\mathrm{S}}}{\mathrm{d}t}
    \end{equation}
    Here, $w_{\mathrm{P}}$ is the wake fraction coefficient. The thrust coefficient curve $K_{\mathrm{T}}(J_{\mathrm{T}})$ can be approximated by the following second order polynomial:
    \begin{equation}
        \label{eq:polynomial_KT}
        K_{\mathrm{T}}(J_{\mathrm{T}}) \approx \kappa_2 J_{\mathrm{T}}^2 + \kappa_1 J_{\mathrm{T}} + \kappa_0  
    \end{equation}
    Next, the authors explain the wave-induced surge force $f(x, t)$. In this paper, the following ITTC1978 wave spectrum is used:
    \begin{equation}
        S_{\mathrm{ITTC}}(\omega) = 172.8\dfrac{H_{1/3}^2}{T_{01}^4 \omega^5} \exp \left( - 691.2 \dfrac{1}{T_{01}^4 \omega^4} \right)
    \end{equation}
    where $\omega$ is the wave frequency, $H_{1/3}$ is the significant wave height, and $T_{01}$ is the mean wave period. The peak frequency $\omega_{\mathrm{WP}}$ is as follows:
    \begin{equation}
        \omega_{\mathrm{WP}}=\left( \frac{4}{5} \frac{691.2}{T_{01}^4} \right)^{\frac{1}{4}}
    \end{equation}
    The wave amplitude at $x_{\mathrm{S}}$ can be calculated with the use of the superposition method as:
    \begin{equation}
        \begin{aligned}
            \zeta_{\omega}(t, x_{\mathrm{S}}(t) = &\displaystyle\sum_{i=1}^{\infty} \sqrt{2S_{\mathrm{W} i}(\omega_{\mathrm{W} i})\mathrm{d}\omega_{\mathrm{W} i}}\\
            &\cdot \cos(\omega_{\mathrm{W} i} t - k_{\mathrm{W} i} x_{\mathrm{S}}(t) + \epsilon_{\mathrm{W} i})
        \end{aligned}
    \end{equation}
    Here, $\omega_{\mathrm{W}i}$: frequency of individual wave, $k_{\mathrm{W} i}$; wave number of individual wave, $\epsilon_{\mathrm{W} i}$: the phase of individual wave. $\epsilon_{\mathrm{W} i} \in [0, 2 \pi)$ is a uniformly distributed random number.

    In our numerical experiments, individual waves with significantly lower or higher frequency components are not required to be generated. Therefore, individual waves inside the frequency range of $\omega_{\mathrm{W}} \in \left[ \frac{1}{2} \omega_{\mathrm{WP}}, 7 \omega_{\mathrm{WP}}\right]$ are used. In this research, the authors assume that the wave-induced surge force can be also calculated by the following superposition method:
    \begin{equation}
        f(x, t) = \displaystyle\sum_{i=1}^{\infty} \mathcal{F}_{\mathrm{W} i} \sin(\omega_{\mathrm{W} i} t - k_{\mathrm{W} i} x_{\mathrm{S}} + \epsilon_{\mathrm{W} i})
    \end{equation}
    Here, $\mathcal{F}_{\mathrm{W}i}$: the amplitude of surge force for individual wave component. Since the $\sin$ function contains the state variable $x_{\mathrm{S}}$, the equation of motion becomes nonlinear as follows:
    \begin{equation}
        \begin{aligned}
            (m + m_x)\frac{\mathrm{d}^2x_{\mathrm{S}}}{\mathrm{d}t^2} = &- R \left(\frac{\mathrm{d}x_{\mathrm{S}}}{\mathrm{d}t} \right) + T \left(\frac{\mathrm{d}x_{\mathrm{S}}}{\mathrm{d}t}, n_{\mathrm{P}} \right)\\
            &+ \displaystyle\sum_{i=1}^{\infty} \mathcal{F}_{\mathrm{W} i} \sin(\omega_{\mathrm{W} i} t - k_{\mathrm{W} i} x_{\mathrm{S}} + \epsilon_{\mathrm{W} i})
        \end{aligned}
    \end{equation}
    Here, it is assumed that the ship resistance and the propeller thrust are balanced with $n_{\mathrm{P}}$ and $\bar{u}_{\mathrm{S}}$.
    \begin{equation}
        - R \left(\bar{u}_{\mathrm{S}}\right) + T \left(\bar{u}_{\mathrm{S}}, n_{\mathrm{P}} \right)=0
    \end{equation}
    The perturbation of surge velocity from $\bar{u}_{\mathrm{S}}$ is defined as $u_{\mathrm{S}}(t)$.
    \begin{equation}
        \frac{\mathrm{d}x_{\mathrm{S}}(t)}{\mathrm{d}t} = \bar{u}_{\mathrm{S}} + u_{\mathrm{S}}(t)
    \end{equation}
    Then, by substituting \Cref{eq:polynomial_resistance} and \Cref{eq:polynomial_KT} into \Cref{eq:EoM_original}, the authors obtain:
    \begin{equation}
        \begin{aligned}
            (m + m_x)\frac{\mathrm{d}^2 x_{\mathrm{S}}(t)}{\mathrm{d}t^2} = &-\sum_{i=1}^{5}\alpha_i u_{\mathrm{S}}^i(t)\\
            &+ \displaystyle\sum_{i=1}^{\infty} \mathcal{F}_{\mathrm{W} i} \sin(\omega_{\mathrm{W}i} t - k_{\mathrm{W}i} x_{\mathrm{S}} + \epsilon_{\mathrm{W}i})
        \end{aligned}
    \end{equation}
    Newly defined symbols in the above equation are as follows:
    \begin{equation}
        \left\{\begin{aligned}
            \alpha_1 = &r_1 + 2 r_2 \bar{u}_{\mathrm{S}} + 3 r_3 \bar{u}_{\mathrm{S}}^2 + 4 r_4 \bar{u}^3_{\mathrm{S}} + 5 r_5 \bar{u}^4_{\mathrm{S}}\\
            &- 2 \bar{u} K_2 (1-t_\mathrm{P})(1-w_\mathrm{P})^2 \rho D_{\mathrm{P}}^2\\
            &- K_1 (1-t_\mathrm{P})(1-w_\mathrm{P}) \rho n D_{\mathrm{P}}^3\\
            \alpha_2 = &r_2 + 3 r_3 \bar{u}_{\mathrm{S}} + 6 r_4 \bar{u}_{\mathrm{S}}^2 + 10 r_5 \bar{u}_{\mathrm{S}}^3\\
            &- K_2 (1-t_\mathrm{P})(1-w_\mathrm{P})^2 \rho D_{\mathrm{P}}^2\\
            \alpha_3 = &r_3 + 4 r_4 \bar{u}_{\mathrm{S}} + 10 r_5 \bar{u}_{\mathrm{S}}^2\\
            \alpha_4 = &r_4 + 5 r_5 \bar{u}_{\mathrm{S}}\\
            \alpha_5 = &r_5
        \end{aligned}\right.
    \end{equation}
    The authors transform to the form of the state equation as follows:
    \begin{equation}
        \left\{\begin{aligned}
            \frac{\mathrm{d}x_{\mathrm{S}}}{\mathrm{d}t} = &\bar{u}_{\mathrm{S}} + u_{\mathrm{S}}(t)\\
            (m+m_x) \frac{\mathrm{d}u_{\mathrm{S}}}{\mathrm{d}t} = &- \sum_{i=1}^{5}\alpha_i u_{\mathrm{S}}^i(t)\\
            &+ \displaystyle\sum_{i} \mathcal{F}_i \sin(\omega_{\mathrm{W}i} t - k_{\mathrm{W}i} x_{\mathrm{S}} + \epsilon_{\mathrm{W} i})
        \end{aligned}\right.
        \label{eq:equation_of_motion_before_approximation}
    \end{equation}

     In the past, analysis for nonlinear stochastic differential equations (SDE: Stochastic Differential Equation) has already been performed in our naval architecture and ocean engineering field. For instance, in the work of Roberts \cite{Roberts1982parametric} and Dostal \cite{dostal2011}, nonlinear stochastic differential equations for parametric rolling are solved by stochastic averaging methods. Recently, Maruyama has analyzed the same problem of parametric rolling using the moment method \cite{Maruyama2023_momentamp}. However, it is difficult to directly apply the stochastic averaging and moment methods to the SDEs treated in this study. This limits the range of theoretical analysis approaches that the authors can use.
    
\section{Numerical examples and their statistical analysis}\label{sec:simulation_results}
\begin{figure*}[htb]
        \centering 
        \includegraphics[width=1.0\hsize]{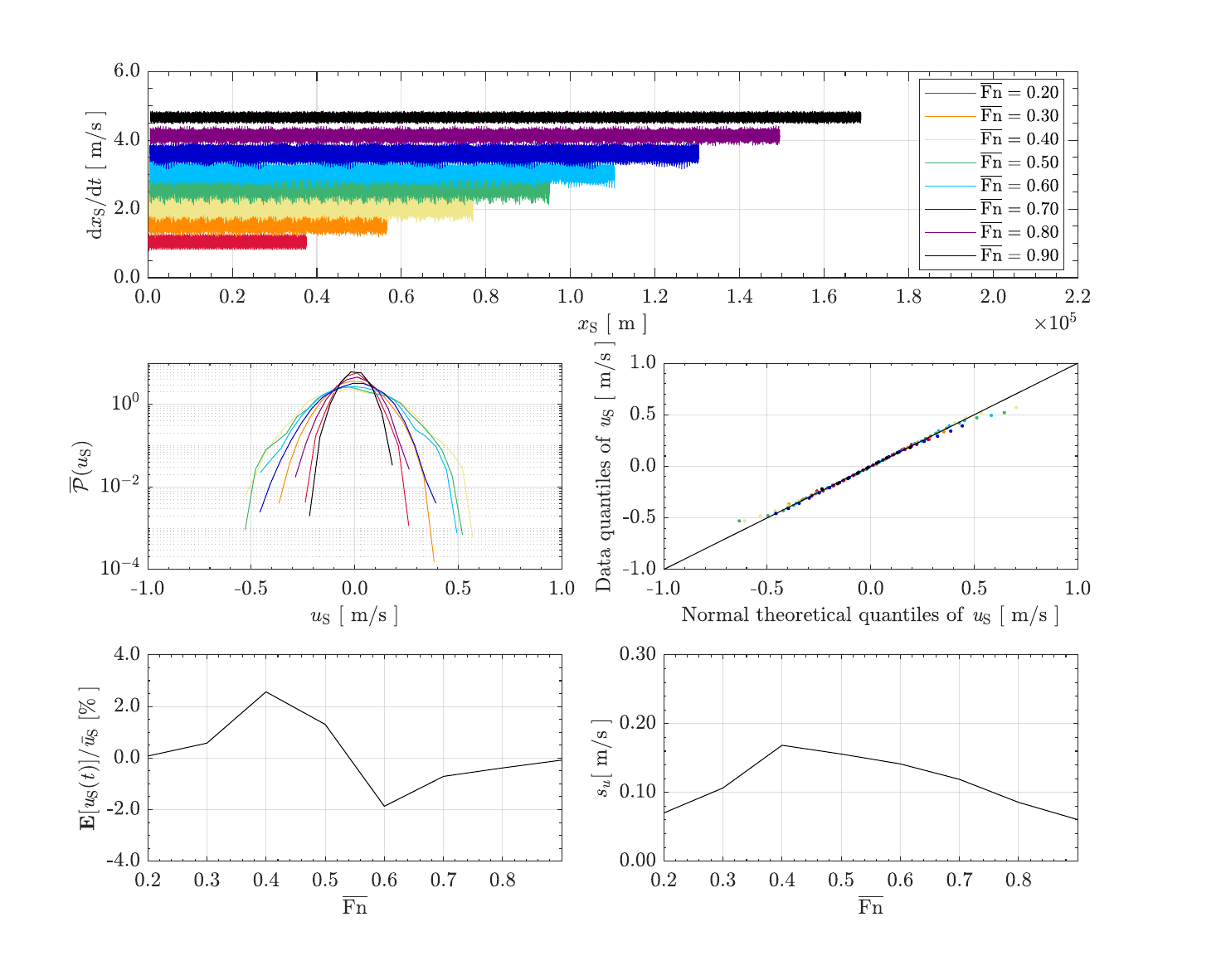}
        \caption{Statistical analysis for several Froude numbers for the equation of motion \Cref{eq:equation_of_motion_before_approximation}} with $T_{01}=1.414~\mathrm{s}$ and $H_{1/3}=0.1~\mathrm{m}$.
        \label{fig:statistical_analysis_result}
    \end{figure*}

    The authors show the numerically solved results of the differential equation derived in \Cref{sec:Equation_of_Motion}. \Cref{fig:statistical_analysis_result} shows the result for 8 different propeller revolution numbers. The uppermost graph shows a phase portrait of the velocity $\mathrm{d}x_{\mathrm{S}}/\mathrm{d}t$ and the traveled distance. The middle left and right graphs represent the PDF and Q-Q plots for the velocity, respectively. The lower left graph shows the percentage of the ratio of the mean velocity in waves to the velocity in calm water, and the lower right graph shows the standard deviation of the velocity. The PDF and Q-Q plots indicate that the velocity can be regarded as almost normally distributed. It can also be seen that there is an increase in average speed of up to $3\%$ when the ship is operated in following seas. It is particularly important to note that the variance of the ship's speed reaches its maximum around $\mathrm{Fn}=0.4$.

    It is noteworthy that the variance of the surge velocity in irregular waves is high for case I of Froude numbers, where surf riding is stable and the convergence to surf-riding is possible in the deterministic case. On the other hand, for case II of Froude numbers, where only surging is possible in the deterministic case, the variance is low.

\section{Approximate solution of the SDE}\label{sec:approximation_solution}
\begin{figure*}[htb]
        \centering 
        \includegraphics[width=1.0\hsize]{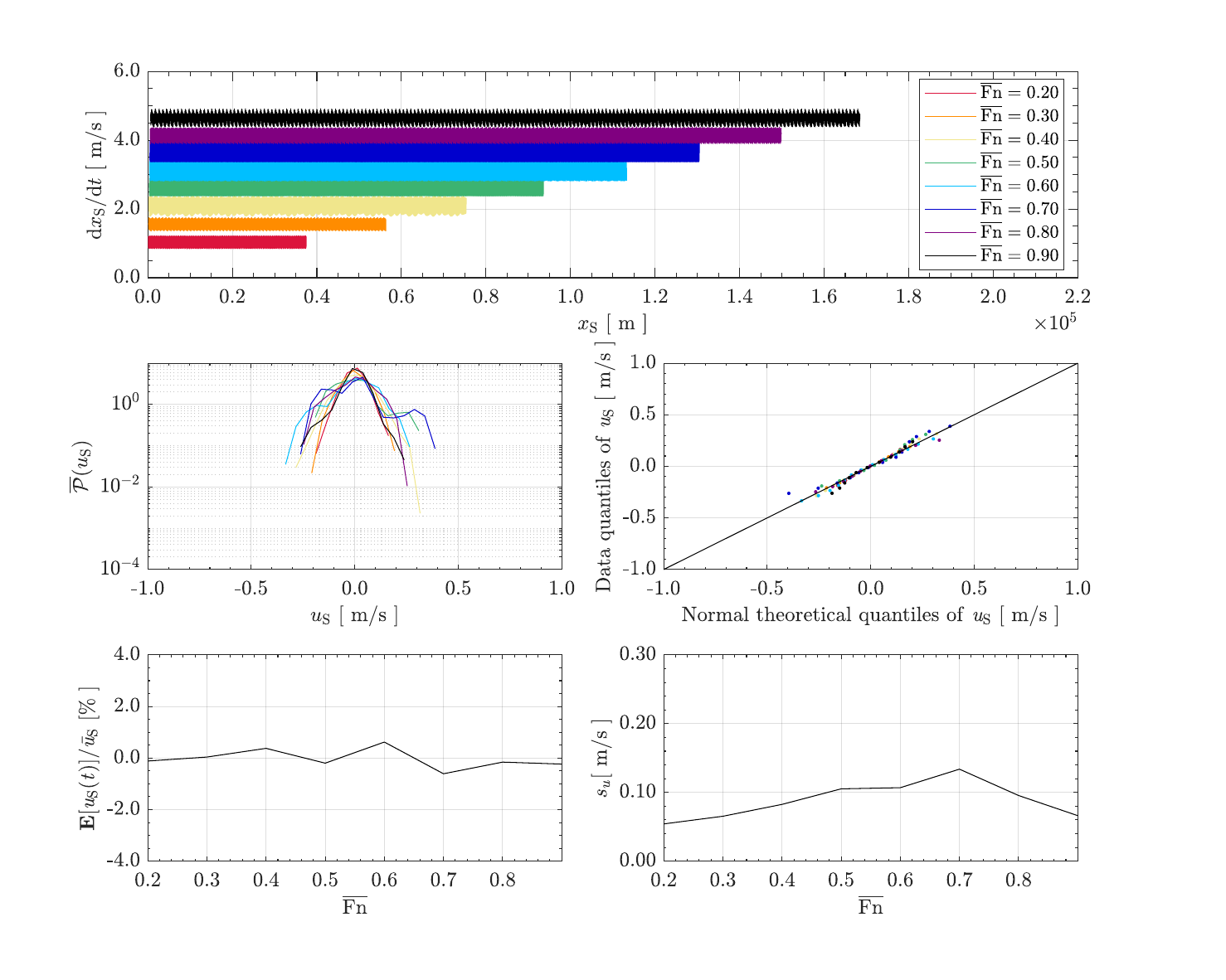}
        \caption{Statistical analysis for several Froude numbers for the approximated equation of motion (\Cref{eq:euqation_of_motion_approximated})}  with $T_{01}=1.414~\mathrm{s}$ and $H_{1/3}=0.1~\mathrm{m}$.
        \label{fig:statistical_analysis_result_for_approximated_system}
    \end{figure*}

    In this section, the authors transform the noise in the equation of motion. 
    \begin{equation}
        \begin{aligned}
            &\displaystyle\sum_{i=1}^{\infty} \mathcal{F}_{\mathrm{W}i} \sin(\omega_{\mathrm{W}i} t - k_{\mathrm{W}i} x_{\mathrm{S}} + \epsilon_{\mathrm{W}i})\\
            =&\displaystyle\sum_{i=1}^{\infty} \mathcal{F}_{\mathrm{W}i} \sin(\omega_{\mathrm{W}i} t + \epsilon_{\mathrm{W}i}) \cos(k_{\mathrm{W}i} x_{\mathrm{S}})\\
            - &\displaystyle\sum_{i=1}^{\infty} \mathcal{F}_{\mathrm{W}i} \cos(\omega_{\mathrm{W}i} t + \epsilon_{\mathrm{W}i}) \sin(k_{\mathrm{W}i} x_{\mathrm{S}})
        \end{aligned}
    \end{equation}
    Now, for the $\sin(k_{\mathrm{W}i} x)$ and $\sin(k_{\mathrm{W}i} x)$ terms, let them be approximated by the following representative values.
    \begin{equation}
        \left\{\begin{aligned}
            \sin(k_{\mathrm{W}i} x_{\mathrm{S}}) &\approx \sin(k_{\mathrm{W}\mathrm{P}} x_{\mathrm{S}})\\
            \cos(k_{\mathrm{W}i} x_{\mathrm{S}}) &\approx \cos(k_{\mathrm{W}\mathrm{P}} x_{\mathrm{S}})
        \end{aligned}\right.
        \label{eq:sinusoidal_comp_approx}
    \end{equation}
    Note that $k_{\mathrm{W}\mathrm{P}}$ is the wave-number corresponding to the peak frequency $\omega_{\mathrm{W} \mathrm{P}}$. This gives the surge force due to the wave as follows
    \begin{equation}
        \begin{gathered}
            \begin{aligned}
                &\displaystyle\sum_{i=1}^{\infty} \mathcal{F}_{\mathrm{W}i} \sin(\omega_{\mathrm{W}i} t - k_{\mathrm{W}i} x_{\mathrm{S}} + \epsilon_{\mathrm{W}i})\\
                \approx &\cos(k_{\mathrm{W}\mathrm{P}} x_{\mathrm{S}}) \mathcal{F}_{\mathrm{W}s}(t)
                -\sin(k_{\mathrm{W}\mathrm{P}} x_{\mathrm{S}}) \mathcal{F}_{\mathrm{W}c}(t)
            \end{aligned}\\
            \text{where} \quad \left\{
            \begin{aligned}
                \mathcal{F}_{\mathrm{W}s}(t) = \displaystyle\sum_{i=1}^{\infty} \mathcal{F}_{\mathrm{W}i} \sin(\omega_{\mathrm{W}i} t + \epsilon_{\mathrm{W}i})\\
                \mathcal{F}_{\mathrm{W}c}(t) = \displaystyle\sum_{i=1}^{\infty} \mathcal{F}_{\mathrm{W}i} \cos(\omega_{\mathrm{W}i} t + \epsilon_{\mathrm{W}i})
            \end{aligned}
            \right.
        \end{gathered}
    \end{equation}
    By utilizing the above, the equation of motion becomes as follows:
    \begin{equation}
        \left\{\begin{aligned}
            \frac{\mathrm{d}x_{\mathrm{S}}}{\mathrm{d}t} = &\bar{u}_{\mathrm{S}} + u_{\mathrm{S}}(t)\\
            (m+m_x) \frac{\mathrm{d}u_{\mathrm{S}}}{\mathrm{d}t} = &- \sum_{i=1}^{5}\alpha_i u_{\mathrm{S}}^i(t)\\
            &+\cos(k_{\mathrm{W}\mathrm{P}} x_{\mathrm{S}}) \mathcal{F}_{\mathrm{W}s}(t)\\
            &-\sin(k_{\mathrm{W}\mathrm{P}} x_{\mathrm{S}}) \mathcal{F}_{\mathrm{W}c}(t)
            \label{eq:euqation_of_motion_approximated}
        \end{aligned}\right.
    \end{equation}
    Similar as in \Cref{fig:statistical_analysis_result}, the obtained equations are numerically solved as in \Cref{fig:statistical_analysis_result_for_approximated_system}.
    From this figure, it can be seen that the change in the variance becomes smaller with the introduction of the approximation {\Cref{eq:sinusoidal_comp_approx}. Such changes in dynamics are basically not negligible. On the other hand, it can be seen at the same time that the qualitative trends before and after the approximation are similar to some extent. Therefore, we will discuss the equations based on this approximation in the following sections.

    Here, as a further approximation, $\mathcal{F}_{\mathrm{W}s}(t)$ and $\mathcal{F}_{\mathrm{W}c}(t)$ are regarded as independent white noise as follows, respectively. This approximation is introduced because colored noise makes analytical investigations more difficult. A similar approach introducing a similar course has been attempted in the previous study on roll motion in transverse waves by Maki et al.~\cite{maki2017Roll}
    
    \begin{equation}
        \left\{\begin{aligned}
            \mathcal{F}_{\mathrm{W}s}(t) \approx \xi_{\mathrm{W}1}(t)\\
            \mathcal{F}_{\mathrm{W}c}(t) \approx \xi_{\mathrm{W}2}(t)
        \end{aligned}\right.
    \end{equation}
    where $\xi_i(t)$ is white noise which has the following properties
    \begin{equation}
        \left\{\begin{aligned}
            \mathbb{E}[\xi_{\mathrm{W}i}(t)] &= 0 \\ \mathbb{E}[\xi_{\mathrm{W}i}(t_1) \xi_{\mathrm{W}i}(t_2)] &= \mathcal{D}_i \delta(t_1 - t_2)
        \end{aligned}\right.
    \end{equation}
    Here, $\mathcal{D}_i$: the intensity of the noise, $\delta(t)$: Dirac's delta function. By using these relationships, the following equation of state is obtained
    \begin{equation}
        \left\{\begin{aligned}
            \frac{\mathrm{d}x_{\mathrm{S}}}{\mathrm{d}t} = &\left(\bar{u}_{\mathrm{S}} + u_{\mathrm{S}}(t) \right)\\
            (m+m_x) \frac{\mathrm{d}u_{\mathrm{S}}}{\mathrm{d}t} = &- \sum_{i=1}^{5}\alpha_i u_{\mathrm{S}}^i(t) \\
            &+ \cos(k_{\mathrm{W}\mathrm{P}} x_{\mathrm{S}}) \xi_{\mathrm{W}1}(t) \\
            &- \sin(k_{\mathrm{W}\mathrm{P}} x_{\mathrm{S}}) \xi_{\mathrm{W}2}(t)
        \end{aligned}\right.
        \label{eq:state_eq_white}
    \end{equation}
    To reformulate this to a stochastic differential equation form, we define the state variable $x(t)$ as:
    \begin{equation}
        x(t) \equiv (x_1(t), x_2(t))^{\top} \equiv (x_{\mathrm{S}}(t), u_{\mathrm{S}}(t))^{\top}
    \end{equation}
    \Cref{eq:state_eq_white} can be modified to Strtonovich-type SDE format as follows:
    \begin{equation}
        \begin{gathered}
            \begin{aligned}
                \mathrm{d}x = &\mu(t, x(t)) \mathrm{d}t + \Sigma(t, x(t)) \circ  \mathrm{d}W(t)
            \end{aligned}\\
            \left\{\begin{aligned}
                \mu(t, x(t)) &\equiv \begin{pmatrix}
                \bar{u}_{\mathrm{S}} + u_{\mathrm{S}}(t) \\
                - \frac{1}{m+m_x} \displaystyle\sum_{i=1}^{5}\alpha_i u^i(t) \\
                \end{pmatrix}
                \\
                \Sigma(t, x(t)) &\equiv \begin{pmatrix}
                0 &0 \\
                \frac{\cos(k_{\mathrm{W}\mathrm{P}} x_{\mathrm{S}})}{m+m_x} &- \frac{\sin(k_{\mathrm{W}\mathrm{P}} x_{\mathrm{S}})}{m+m_x}\\
                \end{pmatrix}\\
                &\equiv \begin{pmatrix}
                \sigma_{11} &\sigma_{12} \\
                \sigma_{21} &\sigma_{22}\\
                \end{pmatrix}\\
                \mathrm{d}W(t) &\equiv \begin{pmatrix}
                \mathrm{d}W_1(t) \\
                \mathrm{d}W_2(t)\\
                \end{pmatrix}
            \end{aligned}\right.
        \end{gathered}
        \label{eq:SDE_eq_white}
    \end{equation}
    Here, $W_1(t)$ and $W_2(t)$ are standard Wiener processes. For further analysis, the authors transform Eq.~\eqref{eq:SDE_eq_white} to an It\^o-type stochastic differential equation as follows:
    \begin{equation}
        \begin{aligned}
            \mathrm{d}x = &\left( \mu(t, x(t)) + \mu_{\mathrm{WZ}}(t, x(t)) \right) \mathrm{d}t\\ 
            &+ \Sigma(t, x(t)) \circ  \mathrm{d}W(t)
        \end{aligned}
        \label{eq:SDE_eq_white_with_WZterm}
    \end{equation}
    In the above equation, $\mu_{\mathrm{WZ}}(t, x(t))$ is the Wong-Zakai's correction term \cite{Wong-Zakai1965}, and it can be calculated as follows:
    \begin{equation}
        \mu_{\mathrm{WZ}i}(t, x(t))= \frac{1}{2} \sum_{j=1}^2 \sum_{k=1}^2 \sigma_{kj} \frac{\partial \sigma_{ij}}{\partial x_k}
    \end{equation}
    Simple manipulations shows that $\mu_{\mathrm{WZ}}(t, x(t))=0$. Therefore hereafter the authors analyze the following SDE:
    \begin{equation}
        \begin{aligned}
            \mathrm{d}x = &\mu(t, x(t)) \mathrm{d}t + \Sigma(t, x(t)) \mathrm{d}W(t)
        \end{aligned}
        \label{eq:SDE_eq_Ito}
    \end{equation}
    The FPK equation for Eq.~\eqref{eq:SDE_eq_Ito} is as follows
    \begin{equation}
        \begin{aligned}
            \frac{\partial}{\partial t}\mathcal{P}(t, x(t)) = -\displaystyle\sum_{i=1}^2 \frac{\partial}{\partial x_i}\mu_i(t,x(t)) \mathcal{P}(t, x(t))\\
            +\frac{1}{2} \displaystyle\sum_{i=1}^2 \displaystyle\sum_{j=1}^2 \frac{\partial^2}{\partial x_i x_j}\left( \displaystyle\sum_{k=1}^2 \sigma_{ik} \sigma_{jk} \mathcal{P}(t,x(t)) \right)
        \end{aligned}
    \end{equation}
    Using the drift and diffusion from Eq.~\eqref{eq:SDE_eq_white} we get the following equation:
    \begin{equation}
        \begin{aligned}
            &\frac{\partial}{\partial t}\mathcal{P}(t, x(t)) =\\
            &-\frac{\partial}{\partial x_{\mathrm{S}}}\left(\left(\bar{u}_{\mathrm{S}} + u_{\mathrm{S}}(t) \right)\mathcal{P}(t, x(t))\right)\\
            &- \frac{\partial}{\partial u_{\mathrm{S}}}\left(-\frac{1}{m+m_x}\left( \sum_{i=1}^{5}\alpha_i u^i(t) \right)\mathcal{P}(t, x(t))\right)\\
            &+\frac{1}{2}\frac{\partial^2}{\partial u_{\mathrm{S}}^2}\left(\frac{\mathcal{D}_1^2 \sin^2(k_{\mathrm{WP}}x_{\mathrm{S}})}{\left(m+m_x\right)^2}+ \frac{\mathcal{D}_2^2 \cos^2(k_{\mathrm{WP}}x_{\mathrm{S}})}{\left(m+m_x\right)^2} \right)\mathcal{P}(t, x(t))
        \end{aligned}
        \label{eq:simpleSDE_eq_Ito}
    \end{equation}
    In this case, we are interested in the stationary probability density function $\bar{\mathcal{P}}(x)$. Therefore, we solve Eq.~\eqref{eq:simpleSDE_eq_Ito} with the left-hand side equated to zero, which yields
    \begin{equation}
        \begin{aligned}
            &0=-\frac{\partial}{\partial x_{\mathrm{S}}}\left(\left(\bar{u}_{\mathrm{S}} + u_{\mathrm{S}} \right) \bar{\mathcal{P}}(x) \right)\\
            &- \frac{\partial}{\partial u_{\mathrm{S}}}\left(-\frac{1}{m+m_x}\left( \sum_{i=1}^{5}\alpha_i u_{\mathrm{S}}^i \right)\bar{\mathcal{P}}(x) \right)\\
            &+\frac{1}{2}\frac{\partial^2}{\partial u_{\mathrm{S}}^2}\left(\frac{\mathcal{D}_1^2 \sin^2(k_{\mathrm{WP}}x_{\mathrm{S}})}{\left(m+m_x\right)^2}+ \frac{\mathcal{D}_2^2 \cos^2(k_{\mathrm{WP}}x_{\mathrm{S}})}{\left(m+m_x\right)^2} \right) \bar{\mathcal{P}}(x)
        \end{aligned}
    \end{equation}
    Hereafter, we proceed with the analysis based on the idea of Mallick and Marcq \cite{mallick2004stochastic}. They assumed that for equations such as Eq.~\eqref{eq:SDE_eq_white}, $\bar{\mathcal{P}}(x)$ is uniformly distributed with respect to $x_{\mathrm{S}}$. Here, the authors show the numerically obtained PDFs of $x_{\mathrm{S}}$ for \Cref{eq:equation_of_motion_before_approximation} and \Cref{eq:euqation_of_motion_approximated} in \Cref{fig:statistical_analysis_result_x} and \Cref{fig:statistical_analysis_result_x_for_approx}, respectively.
    \begin{figure}[htb]
        \centering 
        \includegraphics[width=1.0\hsize]{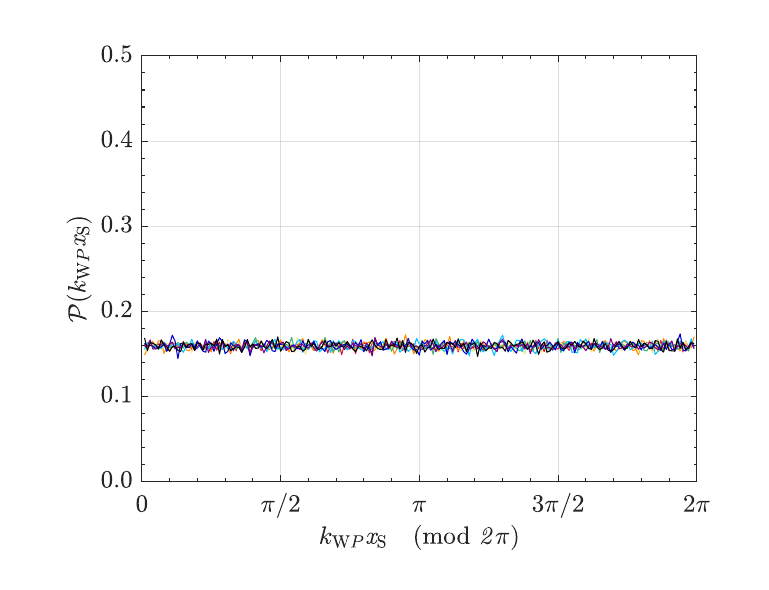}
        \caption{PDF of $k_{\mathrm{WP}}x_{\mathrm{S}}$ for the equation of motion (\Cref{eq:equation_of_motion_before_approximation}})
        \label{fig:statistical_analysis_result_x}
    \end{figure}
    \begin{figure}[htb]
        \centering 
        \includegraphics[width=1.0\hsize]{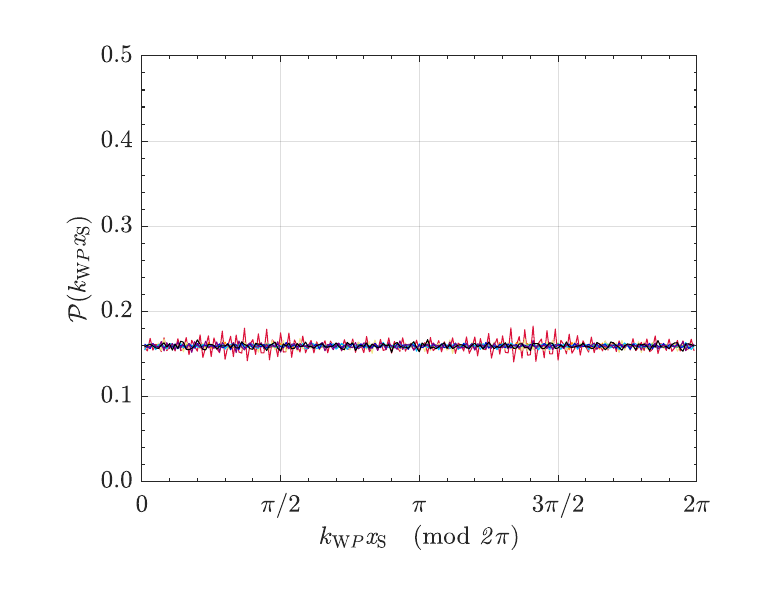}
        \caption{PDF of $k_{\mathrm{WP}}x_{\mathrm{S}}$ for the approximated equation of motion (\Cref{eq:euqation_of_motion_approximated})}
        \label{fig:statistical_analysis_result_x_for_approx}
    \end{figure}
    Thus, we see that the PDF of $x_{\mathrm{S}}$ is uniformly distributed for both equations. Therefore, $\lim_{t \rightarrow +\infty} \mathcal{P}(t,x(t))$ is not a function of position.
    From this, all partial derivatives of $\bar{\mathcal{P}}(x)$ with respect to $x_{\mathrm{S}}(t)$ can be regarded as zero. Therefore, from now on, $\bar{\mathcal{P}}(x)$ will be denoted as $\bar{\mathcal{P}}(u_{\mathrm{S}})$.
    
    Here, the authors average $\sin^2(k_{\mathrm{WP}}x)$ and $\cos^2(k_{\mathrm{WP}}x)$ over one period as:
    \begin{equation}
        \left\{\begin{aligned}
            &\frac{k_{\mathrm{WP}}}{2 \pi} \int_{0}^{\frac{2 \pi}{k_{\mathrm{WP}}}} \sin^2(k_{\mathrm{WP}}x) \mathrm{d}x = \frac{1}{2}\\
            & \frac{k_{\mathrm{WP}}}{2 \pi} 
            \int_{0}^{\frac{2 \pi}{k_{\mathrm{WP}}}} \cos^2(k_{\mathrm{WP}}x) \mathrm{d}x = \frac{1}{2}
        \end{aligned}\right.
    \end{equation}
    This yields the following equation.
    \begin{equation}
        \begin{aligned}
            &\frac{\partial}{\partial u_{\mathrm{S}}}\left(\left( \sum_{i=1}^{5}\alpha_i u_{\mathrm{S}}^i \right)\bar{\mathcal{P}}(u_{\mathrm{S}})\right)\\
            &+ \frac{1}{4}  \frac{\mathcal{D}_1^2 + \mathcal{D}_2^2}{m+m_x} \frac{\partial^2}{\partial u_{\mathrm{S}}^2} \bar{\mathcal{P}}(u_{\mathrm{S}}) = 0
        \end{aligned}
    \end{equation}
    Solving this yields the following probability density function. Note that $\mathcal{D}^2 \equiv \mathcal{D}_1^2 + \mathcal{D}_2^2$.
    \begin{equation}
        \bar{\mathcal{P}}(u_{\mathrm{S}}) = C_{\mathcal{P}}\exp{\left(-\frac{4 \left( m+m_x \right)}{\mathcal{D}^2} \int_{u}\left( \sum_{i=1}^{5}\alpha_i u_{\mathrm{S}}^i \right)\mathrm{d}u_{\mathrm{S}} \right)}
        \label{eq:PDF_us}
    \end{equation} 
    In this PDF $\sum_{i=1}^{5}\alpha_i u_{\mathrm{S}}^i$ is dominated by linear components. Therefore from this, \Cref{eq:PDF_us} has a shape roughly similar to a normal distribution. The theoretical conclusions obtained here are consistent with the results obtained by numerical simulation.

\section{Conclusion}\label{sec:conclusion}
    In this paper, we first organized the statistical aspects of the surf-riding phenomenon. Then, under several approximations, we showed how to calculate the probability density function theoretically. Although the results obtained are based on strong approximations, it was found that the nonlinear back-and-forth oscillations in irregular following waves can be explained from a qualitative point of view. Among the approximations applied, the one with the largest impact is the whitening of the noise. Future work for the authors is to extend the idea of this method to the case of colored noise.

\begin{acknowledgements}
This study was supported by a Grant-in-Aid for Scientific Research from the Japan Society for the Promotion of Science (JSPS KAKENHI Grant \#22H01701). Further, this work was partly supported by the JASNAOE collaborative research program / financial support. The authors are also thankful to Enago (www.enago.jp) for reviewing the English language.
\end{acknowledgements}

%
\section*{Conflict of interest}

    The authors declare that they have no conflict of interest.

\bibliographystyle{spphys}       
\bibliography{main.bib}   

\end{document}